\documentclass[11pt]{amsart}

\usepackage{amssymb,amsmath,amscd,graphicx,
latexsym,amsthm}
\usepackage{amssymb,latexsym,amsmath,amscd,graphicx}
\usepackage[mathscr]{eucal}
  \usepackage[all]{xy}
\def\OO{{\mathcal O}}

\def\F{\mathcal{F}}

\def\G{\mathcal{G}}
\def\H{\mathcal{H}}

\def\cP{\mathcal{P}}
\def\Pic0{{\rm Pic}^0}
\def\Aut0{{\rm Aut}^0}

\def\*{{\underline *}}
\def\utheta{{\underline \theta}}

\def\W{\mathrm{W}}

\def \U{{\mathrm U}}

\theoremstyle{plain}

\newtheorem{theorem}{Theorem}[section]
\newtheorem{theoremalpha}{Theorem}

\newtheorem{proposition/example}[theorem]{Proposition/Example}
\newtheorem{proposition}[theorem]{Proposition}

\newtheorem{corollary}[theorem]{Corollary}

\newtheorem*{theorem*}{Theorem}

\theoremstyle{definition}

\newtheorem{conjecture/question}[theorem]{Conjecture/Question}

\newtheorem{remark/definition}[theorem]{Remark/Definition}
\newtheorem{notation/assumptions}[theorem]{Assumptions/Notation}

\numberwithin{equation}{section}

 \theoremstyle{remark}

\title{ Torsion points on theta divisors and semihomogeneous vector bundles}
 \author[ G. Pareschi]{ Giuseppe Pareschi}

\address{Dipartimento di Matematica,
              Universit\`a di Roma, Tor Vergata\\Italy}
\email{pareschi@mat.uniroma2.it}

 \thanks{The author was partially supported by the MIUR Department of Excellence Project awarded to the Department of Mathematics, University of Rome Tor Vergata, CUP E83C18000100006.}

\begin{document}
\begin{abstract}  We generalize to $n$-torsion  a result of  Kempf's describing $2$-torsion points lying on a theta divisor. This is accomplished by means of certain semihomogeneous vector bundles  introduced and studied by Mukai and Oprea. As an application, we prove a sharp  upper bound  for the number of $n$-torsion points lying on a  theta divisor and show that this is achieved only in the case of products of elliptic curves, settling in the affirmative a conjecture of Auffarth, Pirola and Salvati Manni. 
\end{abstract}
\maketitle

\section{Introduction }

Let $(A,\utheta)$ be a complex $g$-dimensional principally polarized abelian variety. 
 This paper is concerned with the set of $n$-torsion points lying on the theta divisors, where $n$ is any fixed integer $\ge 2$.   
 
 We choose once for all a \emph{symmetric} divisor $\Theta$ representing the polarization. For $x\in A$ we denote  by $t_x:A\rightarrow A$ the translation by $x$, and $\Theta_x$ the effective divisor corresponding to the line bundle $t_x^*\OO_A(\Theta)$ (i.e. $\Theta_x=\Theta-x$). 
    We set
\[\Theta_{x}(n):=\#A[n]\cap \Theta_{x},\]
where $A[n]$ is the group of all $n$-torsion points of $A$.

A   result of Kempf (\cite[Theorem 3]{K1}) asserts that, for $x,y\in A$  the corank of the multiplication map  of global sections
\begin{equation}
\label{m2}H^0(A, t_x^*\OO_A(2\Theta))\otimes H^0(A, t_y^*\OO_A(2\Theta))\longrightarrow H^0(A, t_x^*\OO_A(2\Theta)\otimes t_y^*\OO_A(2\Theta)) 
\end{equation}
coincides with the number $\Theta_{y-x}(2)$ (we refer to \cite[\S2]{PS} for the translation  into the present setting of Kempf's statement, which contains a slight mistake). 
Our first result is an extension  of Kempf's theorem to $n$-torsion points, for arbitrary $n$. This is achieved using certain semihomogeneous vector bundles, denoted $\text{W}_{a,b}$,  introduced and systematically studied by Oprea in \cite{oprea} (as a consequence of Mukai's theory of semihomogeneous vector bundles, \cite{mukai-semi}). When  $a$ and $b$ are coprime positive integers, the  $\W_{a,b}$'s are defined as simple, semihomogeneous and symmetric vector bundles  such that 
\begin{equation}\label{defining}
 \mathrm{rk}\W_{a,b}=a^g\qquad \hbox{and}\qquad\det\W_{a,b}=\OO_A(\Theta)^{a^{g-1}b}.
\end{equation}
If $a$ is odd there is a unique such vector bundle, while if $a$ is even  they are not unique when $g\ge 2$. We refer to the next section for generalities on such vector bundles. Our generalization of Kempf's theorem (which is recovered for $a=b=1$) is the following

\begin{theoremalpha} \label{kempf-gen}
 Let $a$ and $b$ be coprime positive integers.   Let $\W_{a,a+b}$ and $\W_{b,a+b}$ be two vector bundles as above. For $x,y\in A$  the number $\Theta_{y-x}(a+b)$ is equal to the corank of the multiplication map of global sections
\begin{equation}\label{W}m_{a,b}(x,y): H^0(A,t_x^*\W_{a,a+b})\otimes H^0(A, t_y^*\W_{b,a+b})\longrightarrow H^0(A,t_x^*\W_{a,a+b}\otimes t^*_y\W_{b,a+b})
\end{equation}
\end{theoremalpha}

 (As  it is easy to check,  the source and target of the above map have the same dimension, namely $(a+b)^{2g}$.) Note that  if $a$ or $b$ are even, say $a$, the above map $m_{a,b}$ depends on a choice of a vector bundle $\W_{a,a+b}$, but we will neglect this in the notation. A special role will be played by
the particular case
\begin{equation}\label{n-1}
m_{1,n-1}(x,y): H^0(A,t_x^*\OO_A(n\Theta))\otimes H^0(A, t_y^*\W_{n-1,n})\longrightarrow H^0(A,t_x^*\OO_A(n\Theta)\otimes t^*_y\W_{n-1,n})
\end{equation}
obtained for $a=1$ and $n:=a+b$.

In view of Theorem \ref{kempf-gen}, it is useful to consider criteria for the surjectivity of the multiplication of global sections of semihomogeneous vector bundles, in analogy with well known classical theorems   for line bundles (due to Mumford, Koizumi, Kempf and others, see e.g. \cite[\S 7.2]{BL} \cite[\S 6.2]{Kempf} , \cite[\S 8]{JP}). In fact an optimal result in this direction  was already proved years ago by Popa and the author  (\cite[Theorem 7.30]{PP2}). We restate  it more expressively as Theorem \ref{pp} below.
In turn this is an ingredient of the proof of the following lower bound for the rank of the multiplication maps appearing in Theorem \ref{kempf-gen}

\begin{theoremalpha}\label{cor}  
 In the notation of Theorem \ref{kempf-gen}, \ $\text{rk}\, (m_{a,b}(x,y))\ge((a+b)^2-1)^g$  for all $x,y\in A$ .
\end{theoremalpha}

The last result of this note, in fact our original motivation, is the proof of the following  conjecture of Auffarth, Pirola and Salvati Manni  on the maximal  number of $n$-torsion points on a theta divisor  (see \cite{APS}).  The case $n=2$, which was conjectured earlier by Marcucci and Pirola (\cite{MP}), was  proved by Salvati Manni and the author in \cite{PS}. 
\begin{theoremalpha} \label{main}   For all $x\in A$
$$\Theta_{x}(n)\leq n^{2g} -(n^2-1)^g.$$
Moreover equality  holds  if and  only if  $A$  is a product of  elliptic  curves and \ $\OO_A(\Theta_{x})=\boxtimes_{i}\OO_{E_i}(z_i)$ \ where $z_i$ are  n-torsion points.
\end{theoremalpha}

Note that the upper bound of the statement  is just the combination of Theorems \ref{kempf-gen} and  \ref{cor}. The remaining part is proved in Section \ref{proof}. 

 It would be interesting to have effective results along these lines for irreducible principal polarizations. To this purpose it should be kept in mind that, thanks to a recent result of Auffarth and Codogni (\cite{AC}), there are irreducible theta divisors containing abelian subvarieties   of dimension  up to $g-2$, hence containing at least $n^{2(g-2)}$ \ $n$-torsion points for all $n$.  On the other hand, by Raynaud's theorem on the Manin-Mumford conjecture,  the overall number of torsion points contained in a theta divisor is finite unless it contains translates of positive-dimensional abelian subvarieties by torsion points. 

\noindent\textbf{Acklowledgements. } The author thanks:  Dragos Oprea for pointing out a gap in an earlier draft of this paper, Mihnea Popa for many conversations about semihomogeneous vector bundles a long time ago, Riccardo Salvati Manni for his encouragement and many discussions and suggestions, and the referee for very accurate remarks.


\section{Preliminaries on the vector bundles $\W_{a,b}$}\label{prelim} 

Here we recall some basic facts about  the vector bundles $\W_{a,b}$  introduced by Oprea in \cite{oprea}. Let $(A,\utheta)$ be a $g$-dimensional p.p.a.v. and  let $\Theta$ be a fixed symmetric theta divisor. For a pair of coprime positive integers $a$ and $b$ we consider  simple semihomogeneous vector bundles $\W$ such that
\begin{equation}\label{defining1}
 rk\W=a^g\qquad \det\W=\OO_A(\Theta)^{a^{g-1}b} \>.
\end{equation}
   Vector bundles with the above properties  do  exist and they are unique up to tensorization with  an $a^g$-torsion line bundle (\cite[Theorem 7.11 and Remark 7.13]{mukai-semi}).  Denoting by $a_A:A\rightarrow A$ the isogeny $x\mapsto ax$, we have that
\begin{equation}\label{fund}
a_A^*\W\cong \bigl(\OO_A(\Theta)^{ab}\bigr)^{\oplus a^g}
\end{equation}
(see \cite[2.3.1]{oprea}). Moreover such $\W$'s satisfy \emph{the index theorem with index $0$} (IT(0) for short), meaning that $h^i(\W_{a,b}\otimes \alpha)=0$ for all $i>0$ and $\alpha\in \widehat A$ (we denote $\widehat A=\Pic0 A$ the dual abelian variety).  Recalling that the degree of the isogeny $a_A$ is $a^{2g}$, it follows that
\begin{equation}\label{rk-etc}
 h^0(A,\W)=\chi(W)=a^g\bigl(\frac{b}{a}\bigr)^g=b^g\>.
\end{equation}
Another useful fact about  the vector bundles satisfying the above condition (\ref{defining1}) is that they are globally generated as soon as $b>a$. This follows from the criterion asserting  that a vector bundle $E$ on $A$ is globally generated as soon as $E(-\Theta)$ is IT(0) (\cite[Theorem 2.1]{pareschi}). Indeed  (\ref{fund})  yields that $\W(-\Theta)$ is IT(0) if and only if $b>a$.

 For odd $a$, imposing the supplementary condition  that $\W$ is \emph{symmetric}, i.e. $(-1)_A^*\W\cong \W$, it turns out that there is  a \emph{unique} such vector bundle up to isomorphism (\cite[\S2.1]{oprea}). It is denoted $\W_{a,b}$.
Also for even $a$  such  symmetric vector bundles do exist (for example the dual of the Fourier-Mukai transform of the vector bundle $\W_{b,a}$) but they are  not unique for $g\ge 2$. For even $a$ Oprea defines a unique such 
vector bundle $\W_{a,b}$ by means of the Schr\"odinger representation (\cite[\S2.1 and (16)]{oprea}). However this is not  important for our purposes, since we will consider \emph{any} simple symmetric semihomogeneous vector bundle $\W$ satisfying (\ref{defining1}).  We  denote $\mathcal{W}_{a,b}$ the set of all isomorphism classes of such simple symmetric semihomogeneous vector bundles and a
vector bundle $\W\in\mathcal W_{a,b}$ will be usually denoted $\W_{a,b}$.\footnote{Here our notation  differs from the one of Oprea, as he denotes $\W_{a,b}$ the unique vector bundle in $ \mathcal{W}_{a,b}$ defined, as mentioned above, via the Schr\"odinger representation.} 

 We consider the subgroup 
\[\Sigma(\W_{a,b})=\{\alpha\in\widehat A\>|\> \W_{a,b}\otimes \alpha\cong\W_{a,b}\}.\]
We have that, independently on the parity of $a$, 
\begin{equation}\label{subgroup}\Sigma(\W_{a,b})=\widehat A[a]
\end{equation}
 ($a$-torsion line bundles, see \cite[Corollary 7.2]{mukai-semi}). 
 
 Given $\W_{a,b}\in \mathcal W_{a,b}$, the other vector bundles, say $\W_{a,b}^\prime$ (possibly isomorphic to $\W_{a,b}$) whose isomorphism class lies in $\mathcal{W}_{a,b}$ are those of the form
\[\W_{a,b}^\prime\cong \W_{a,b}\otimes \beta\]
for $\beta\in \widehat A[a^g]$ such that $(-1)_A^*(\W_{a,b}\otimes\beta)\cong  \W_{a,b}\otimes\beta^{-1}\cong\W_{a,b}\otimes\beta$. Therefore
 $\beta^2\in \widehat A[a]$. This, together with the condition $\beta\in \widehat A[a^g]$ implies that, if $a$ is odd (or $g=1$) then $\beta\in\widehat A[a]=\Sigma(\W_{a,b})$. Hence, as mentioned above, there is a unique such an isomorphism class (\cite[\S2.1]{oprea}). 

 In the proof of Theorem \ref{kempf-gen} the following (slight variant of a) result of Oprea will be in use.  
 For $a$ and $b$ coprime positive integers
one considers the isogeny 
\[\mu_{b, a}:A\times A\rightarrow A\times A\quad (z,t)\mapsto (bz+at, z-t).\]

\begin{proposition}\label{maledetta} \emph{(Oprea)} 
Keeping the above notation and assumptions, 
 given a pair of vector bundles $(\W_{a,a+b},\W_{b,a+b}) \in\mathcal{W}_{a,a+b}\times\mathcal W_{b,a+b}$ there is a vector bundle    $W_{ab,1}\in\mathcal W_{ab,1}$ such that   
\begin{equation}\label{wirt2} 
\mu_{b,a}^*(\W_{ab,1}\boxtimes \OO_A(\Theta))\cong  \W_{a,a+b}\boxtimes \W_{b,a+b}.
\end{equation}
\end{proposition}

\proof 
Let $a$ and $b$ be coprime positive integers. For $a$ and $b$ both odd (or $g=1$ and arbitrary $a$ and $b$, see \cite{a})  all vector bundles appearing in the statement are unique and the Proposition is exactly  Oprea's \cite[Proposition 1]{oprea}. 
Next, we assume that  $g\ge 2$ and $a$ and $b$  are still coprime, but one of them, say $a$,  is even.   
We fix a  reference bundle $\overline\W_{ab,1}\in \mathcal{W}_{ab,1}$.   Oprea's argument  still  proves that  the determinant of  $\mu_{b,a}^*\bigl(\overline\W_{ab,1}\boxtimes \OO_A(\Theta)\bigr)$  is equal to the one of $\W_{a,a+b}\boxtimes\W_{b,a+b}$, and that
\begin{equation}\label{wirt3} 
\mu_{b,a}^*(\overline\W_{ab,1}\boxtimes \OO_A(\Theta)) \cong (\W_{a,a+b}\otimes \delta)\boxtimes (\W_{b,a+b} \otimes \gamma)
\end{equation}
for suitable  $(\delta,\gamma)\in \widehat A[a^g]\times\widehat A[b^g]$. We claim that, moreover, both the vector bundles $E:=\W_{a,a+b}\otimes \delta$ and $F:=\W_{b,a+b}\otimes \gamma$ are symmetric. Indeed  an immediate computation shows that $(-1_A,1_A)\circ\mu_{b,a}= (1_A,-1_A)\circ \mu_{b,a}\circ (-1_A,-1_A)$. Therefore, since $\overline\W_{ab,1}$ and $\OO_A(\Theta)$ are both symmetric, pulling back $\overline\W_{ab,1}\boxtimes \OO_A(\Theta)$ under the morphism $(-1_A,1_A)\circ\mu_{b,a}$ we get that $(-1_A,1_A)^*(E\boxtimes F)\cong (1_A,-1_A)^*(E\boxtimes F)$. This proves what was claimed. Hence, since $b$ is odd,  $\W_{b,a+b}\cong \W_{b,a+b}\otimes\gamma$, i.e. $\gamma\in\widehat A[b]$. Moreover $\delta^2\in \Sigma(\W_{a,a+b})=\widehat A[a]$. 

To conclude the proof, we show that we can replace $\overline \W_{ab,1}$ with another vector bundle in $\W_{ab,1}\in \mathcal W_{ab,1}$ such that (\ref{wirt2}) is satisfied.  For any $\alpha\in \widehat A$ we have that $\mu_{b,a}^*(\alpha\boxtimes \OO_A)=(\alpha^b,\alpha^a)$. Therefore  
\[\mu_{b,a}^*\bigl((\overline\W_{ab,1}\otimes\alpha)\boxtimes \OO_A(\Theta)\bigr) \cong (\W_{a,a+b}\otimes \delta\otimes \alpha^b)\boxtimes (\W_{b,a+b} \otimes \gamma\otimes\alpha^a).\]
Taking any $\alpha$ such that $\alpha^b=\delta^{-1}$ we have that $\alpha^2\in\widehat A[ab]$ (hence also  $\alpha\in\widehat A[(ab)^g]$ as soon as $g> 1$). Therefore  $\overline\W_{ab,1}\otimes\alpha\in\mathcal{W}_{ab,1}$. As above, by uniqueness when $b$ is odd, we have that $\W_{b,a+b}\cong \W_{b,a+b}\otimes \gamma\otimes \alpha^a$.  Hence  
\[
\mu_{b,a}^*\bigl((\overline\W_{ab,1}\otimes\alpha)\boxtimes \OO_A(\Theta)\bigr) \cong \W_{a,a+b}\boxtimes \W_{b,a+b}.
\]

\section{Proof of Theorem \ref{kempf-gen} }\label{kempf-section} 

We essentially follow Kempf's argument, with some simplifications. Let
 \[ (\W_{a,a+b}, \W_{b,a+b})\in \mathcal W_{a,a+b}\times\mathcal W_{b,a+b}.\]
  To render the argument more transparent we first prove the result for $(x,y)=(0,0)$. 

The multiplication map $m_{a,b}(0,0)$ is the map $H^0(r_\Delta)$, where $r_\Delta$ is  the restriction to the diagonal
\[ 
r_\Delta: \W_{a,a+b}\boxtimes \W_{b,a+b}\rightarrow \bigl(\W_{a,a+b}\boxtimes \W_{b,a+b}\bigr)_{|\Delta}.
\]
We apply Proposition \ref{maledetta}, ensuring that
\[
\mu_{b,a}^*\bigl((\W_{ab,1})\boxtimes \OO_A(\Theta)\bigr)\cong  \W_{a,a+b}\boxtimes \W_{b,a+b}.
\]
Moreover, since $\mu_{b,a}^{*}(\OO_{A\times \{0\}})=\OO_\Delta$, it follows that $r_\Delta=\mu_{b,a}^*(\rho)$, where $\rho$ is the restriction map
\[\rho: \W_{ab,1}\boxtimes \OO_A(\Theta)\rightarrow \bigl(\W_{ab,1}\boxtimes \OO_A(\Theta)\bigr)_{|A\times\{0\}}.\
\]
(By the way we note that, since the restricted map $(\mu_{b,a})_{|\Delta}:\Delta\rightarrow A\times\{0\}$ is identified to  the isogeny $(a+b)_A:A\rightarrow A$, it follows that $\W_{a,a+b}\otimes \W_{b,a+b}\cong (a+b)_A^*\W_{ab,1}$). 

The kernel of the isogeny $\mu_{b,a}$ is $\Delta[a+b]:=\{(z,z)\>|\>(a+b)z=0\}$. Therefore, since 
\[
r_\Delta=\mu_{b,a}^*(\rho),
\]
 the multiplication map $m_{a,b}(0,0)=H^0(r_\Delta)$ decomposes as 
\[\bigoplus_{\alpha\in\widehat A[a+b]}\Bigl(H^0(A\times A, (\W_{ab,1}\otimes P_\alpha)\boxtimes (\OO_A(\Theta)\otimes P_\alpha))\rightarrow H^0(A\times A, ((\W_{ab,1}\otimes P_\alpha)\boxtimes (\OO_A(\Theta)\otimes P_\alpha))_{|A\times\{0\}}\Bigr)\>.\]
Via the isomorphism induced by the principal polarization $A\rightarrow \widehat A$, 
the above can be written as 
\[\bigoplus_{z\in A[a+b]}\Bigl(H^0(A\times A, t_z^*\W_{ab,1}\boxtimes  t_z^*\OO_A(\Theta))\buildrel{\lambda_z}\over\longrightarrow H^0(A\times A, (t_z^*\W_{ab,1}\boxtimes t_z^*\OO_A(\Theta))_{|A\times\{0\}})\Bigr)\>.\]
Notice that, by (\ref{rk-etc}), $H^0(A,\W_{ab,1})=1$, so that the individual maps $\lambda_z$ appearing above are maps of $1$-dimensional vector spaces. Hence the assertion of the theorem follows from the fact that  the scalar $\lambda_z$ vanishes if and only if  $A\times\{0\}\subset A\times \Theta_{z}$,  i.e. $z\in\Theta$.

In the general case the proof is similar. In the first place, applying $t_{-x}^*$ we can assume that $x=0$. Via translation on the second factor we identify the  map $m_{a,b}(0,y)$ of the statement to the map 
\begin{equation}\label{variant} H^0(A,\W_{a,a+b})\otimes H^0(A, \W_{b,a+b})\longrightarrow H^0(A,\W_{a,a+b}\otimes t^*_y\W_{b,a+b}).
\end{equation}
This  is the $H^0$ of the restriction map 
\begin{equation}\label{variant2}
r_{\Delta_y}: \W_{a,a+b}\boxtimes \W_{b,a+b}\rightarrow \bigl(\W_{a,a+b}\boxtimes \W_{b,a+b}\bigr)_{|\Delta_y}
\end{equation}
where $\Delta_y=d^{-1}(y)$ (here $d$ is the difference map $A\times A\rightarrow A$, $(z,t)\mapsto z-t$). We have that $\OO_{\Delta_y}=\mu_{b,a}^{*}(\OO_{A\times\{y\}} )$. The rest of the proof is unchanged.

\section{Multiplication of global sections of semihomogeneous vector bundles}\label{multip-section}

\noindent\textbf{A surjectivity criterion for multiplication maps. }  We recall   \cite[Theorem 7.30]{PP2}, mentioned in the Introduction.  In the case of interest for this paper, namely semihomogeneous vector bundles whose first Chern class is a power of a principal polarization,   it can be  stated as follows. 
Following Mukai (\cite{mukai-semi}), for a vector bundle $E$ on $A$ we write
\begin{equation}\label{delta}
\delta_E=\frac{c_1(E)}{\text{rk}\,(E)}\in NS(A)\otimes \mathbb Q \>.
\end{equation} 
If $c_1(E)$ is a multiple of $\utheta$ we denote also $\mu_E$ the rational number defined by
\[\delta_E=\mu_E\utheta \>.
\]

\begin{theorem}\label{pp}(Pareschi-Popa) Let $E$ and $F$ be semihomogeneous vector bundles on $A$ such that $c_1(E)$ and $c_1(F)$ are multiples of $\utheta$. If  \[\mu_F >1\qquad \hbox{and}\qquad
\mu_E>\frac{\mu_F}{\mu_F-1}
\]
then the multiplication map of global sections
\[
H^0(A,E)\otimes H^0(A,F)\rightarrow H^0(A,E\otimes F)
\]
is surjective. 
\end{theorem}
Note that for line bundles one recovers the classical fact that the multiplication map of a second power and a third power of a line bundle representing $\utheta$ is surjective.

Here we show that Theorem \ref{pp} is just the restatement of  \cite[Theorem 7.30]{PP2},  asserting that the multiplication map as in the statement  is surjective as soon as both $E(-\Theta)$ and $F(-\Theta)$ satisfy IT(0) and  
\begin{equation}\label{mascherata}
\delta_{E(-\Theta)}+\delta_{\widehat{\Phi}_\cP(F(-\Theta))}>0
\end{equation}
where $\widehat{\Phi}_{\cP}: D(A)\rightarrow D(\widehat A)$ is the Fourier-Mukai transform associated to the Poincar\'e bundle.  Let us explain how to get  the statement of Theorem \ref{pp} from this.   In the first place we recall that, for a semihomogeneous vector bundle $G$, the IT(0) condition is equivalent to $\delta_G>0$. This follows, for example, from \cite[Lemma 6.11]{mukai-semi}, stating that 
\begin{equation}\label{homog}
r_A^*G\cong (\det G)^r\otimes H
\end{equation}
 where   $r:=\text{rk}\,G$, $r_A$ denotes, as usual,  the isogeny $x\mapsto rx$, and $H$ is a homogeneous vector bundle (indeed a homogeneous vector bundle is a direct sum of vector bundles of the form $\U\otimes \alpha$, where $\alpha\in \widehat A$ and $\U$ is a \emph{unipotent} vector bundle, namely a vector bundle having a filtration  $0=\U_0\subset \U_1\subset\cdots\subset\U_{n-1}\subset \U_n=\U$, with $\U_i/U_{i-1}\cong\OO_A$ for $i=1,\dots ,n$, see \cite[Theorem 4.17]{mukai-semi}).

 Using the formulas 
    \begin{equation}\label{sum}
   \delta_{F\otimes G}=\delta_F+\delta_G\qquad\hbox{ and }\qquad \delta_{G^\vee}=-\delta_G,
   \end{equation}
   it follows that the condition that $\mu_F-1>0$, i.e. the first hypothesis of Theorem \ref{pp}, is equivalent to the fact that $F(-\Theta)$ satisfies IT(0).
    
If this is the case then, by base change, the complex $\widehat{\Phi}_\cP(F(-\Theta))$ is a sheaf in cohomological degree $0$, in fact a locally free sheaf. 
 Next, we recall that, for a semihomogenous vector bundle $G$ satisfying IT(0), the vector bundle $\widehat{\Phi}_\cP(G)$ is again semihomogeneous (this follows  from the fact that $\widehat{\Phi}_\cP$ exchanges translation with tensorization with a line bundle in $\widehat A$, see \cite[(3.1)]{mukai}). 
 Finally we claim that if $G$ is such that $c_1(G)$ is a multiple of $\utheta$ then also $c_1(\widehat{\Phi}_\cP(G))$ is a multiple of $\utheta$ and  the following beautiful formula holds:
 \begin{equation}\label{nice}
 \mu_{\widehat{\Phi}_\cP(G)}=-\frac{1}{\mu_G}.
 \end{equation}
 This translates the hypothesis (\ref{mascherata}) into the numerical condition
 \[\mu_E-1-\frac{1}{\mu_F-1}>0,\]
 i.e. the second inequality in the hypothesis of Theorem \ref{pp}.

  Finally, we briefly indicate the proof of (\ref{nice}).  This is certainly well known to the experts but we couldn't find an explicit reference. 
We recall that, for $\lambda\in\mathbb Q$,   $r_A^*(\lambda\utheta)=r^2\lambda\utheta$. Therefore from (\ref{homog})  it follows that $ch(G)=r\exp (\mu_G\utheta)$. Then a well known calculation using GRR and the Fourier-Mukai transform at the level of Chow rings modulo numerical equivalence (see e.g. the proof of \cite[Lemma 2]{oprea}) shows that $ch(\widehat{\Phi}_\cP(G))=r(\mu_G)^{g}\exp (-\mu_G^{-1}\utheta)$. Therefore 
 $\delta_{\widehat{\Phi}_\cP(G)}=-\mu_G^{-1}\utheta$.

\noindent\textbf{Proof of Theorem \ref{cor}. }  By Theorem \ref{kempf-gen} for fixed $x,y\in A$ the maps $m_{a,b}(x,y)$ have the same rank  for all $a,b$ (coprime) with  $a+b=n$ (and for all representatives $(\W_{a,a+b}, \W_{b,a+b})\in\mathcal {W}_{a,a+b}\times\mathcal{W}_{b,a+b}$).  Therefore it is enough to prove the statement for $(a,b)= (n-1,1)$. As in the proof of Theorem \ref{kempf-gen}, we can furthermore assume that $x=0$. Let us fix $y\in A$. For general $z\in A$ we consider the commutative diagram
 \begin{equation}\label{diagram}\xymatrix{H^0(\W_{n-1,n})\otimes H^0(t_y^*\OO_A(n\Theta))\otimes H^0(t_z^*\W_{n-1,n})\ar[rr]\ar[dd]^{m_{n-1,1}(0,y)\otimes \,\mathrm{id}}&&H^0(\W_{n-1,n})\otimes H^0(t_y^*\OO_A(n\Theta)\otimes t_z^*\W_{n-1,n})\ar[dd]\\ \\
  V_n(y)\otimes H^0(t_z^*\W_{n-1,n})\ar[rr]&&H^0(\W_{n-1,n}\otimes t_y^*\OO_A(n\Theta)\otimes t_z^*\W_{n-1,n})}
  \end{equation}
  where $V_n(y)$ denotes the image of the map $m_{n-1,1}(0,y)$. 
  By Theorem \ref{kempf-gen} the top horizontal map is surjective for general $z\in A$. By Theorem \ref{pp} the right vertical map is surjective for all $z\in A$ (here we use that $\mu_{\W_{a,b}}=\frac{b}{a}$ and (\ref{sum})). Therefore the bottom horizontal map is surjective for general $z\in A$. Hence 
  \[
  \dim V_n(y)\ge \frac{\chi(\W_{n-1,n}^{\otimes 2}\otimes\OO_A(n\Theta))}{\chi(\W_{n-1,n})} \>.
  \]
  By (\ref{rk-etc}) we have that $\chi(\W_{n-1,n})=n^g$. Using  (\ref{fund}) one gets easily that $\chi(\W_{n-1,n}^{\otimes 2}\otimes \OO_A(n\Theta))=(n-1)^gn^g(n+1)^g$. The result follows.\qed


\section{Proof of Theorem \ref{main}}\label{proof}

 It is easy to check that the bound of Theorem \ref{main} is attained    by   line bundles of the form 
  $ \OO_A(\Theta_x)= \boxtimes_{i}\OO_{E_i}(z_i)$ on products of elliptic curves $E_i$,  
   where the $z_i$'s are $n$-torsion points. Conversely, the second part  of Theorem \ref{main} asserts that this is the only case. To prove this, the main point consists in showing that if the bound is attained then $\Theta$ must be reducible and therefore the p.p.a.v. $(A,\utheta)$ decomposes as a product of lower dimensional p.p.a.v's.

\noindent\textbf{Proof for $\mathbf{g>2}$. }     According to Theorem \ref{kempf-gen}, what we need to show is that in the irreducible case the rank of the multiplication maps $m_{n-1,n}(0,y)$ is $>(n^2-1)^2$ for all $y$. 
   
 For $y\in A$ and $n>1$ let us consider the following divisor
\begin{equation}\label{E}
E_{y,n}:=\sum_{\eta\in A[n]}\Theta_{y+\eta}.
\end{equation}
An immediate consequence of Theorem \ref{kempf-gen} is the following (where, for reasons apparent in what follows, the notation $(x,y)$ in the statement of Theorem \ref{kempf-gen} has been switched to $(y,z)$):
\begin{corollary}\label{aggiunto} Let $y\in A$. The map 
\[
m_{a,b}(y,z): H^0(A,t_y^*W_{a,a+b})\otimes H^0(A, t_z^*\W_{b,a+b})\longrightarrow H^0(A,t_y^*\W_{a,a+b}\otimes t^*_z\W_{b,a+b})
\]
is singular if and only if $z\in \text{Supp}\,E_{y,a+b}$.
\end{corollary}

In order to prove Theorem \ref{main} we see the  map of Corollary \ref{aggiunto} and all  maps of diagram  (\ref{diagram}) as the fiberwise maps of maps of locally free sheaves. This is well known and it is done as follows.
 Following  \cite{pareschi} (see also \cite{PP1}), given two coherent sheaves $\F$ and $\G$, we define their (derived) skew Pontryagin product
\[ 
\F\mathbin{\hat {*}}\G:=d_*(\F\boxtimes \G)
\]
where $d:A\times A\rightarrow A$ is the difference map. We  make the simplifying assumption that both  sheaves $\F$ and $\G$ are   locally free, semihomogeneous, and they satisfy IT(0) (these conditions will be always satisfied by the sheaves appearing in what follows).  By \cite[Proposition 2.9]{PP1}  the IT(0) condition for $\F$ and $\G$ implies that also the vector bundle $\F\otimes \G$ satisfies IT(0). In turn, this implies by base change  that:   $R^id_*(\F\boxtimes \G)=0$ for $i\ne 0$,  $\F\hat *\G$ is a locally free sheaf (in degree  0), and
\[ d_*(F\boxtimes G)\otimes \mathbb C(z)\cong H^0(A\times A, (F\boxtimes G)_{|\Delta_z})=H^0(A,\F\otimes t_z^*\G)\]
for all $z\in A$ (see (\ref{variant2})). Thus the multiplication map of global sections 
 \[
H^0(A, \F)\otimes H^0(A, t_z^*\G)\rightarrow H^0(A, F\otimes t_z^*G),
 \]
 is naturally identified, via the isomorphism 
 \[
 \mathrm{id}\otimes t_z^*:  H^0(A, \F)\otimes H^0(A, \G)\rightarrow H^0(A, \F)\otimes H^0(A, t_z^*\G),
 \] to the fiber map at $z$ of the map of $\OO_A$-modules:
\begin{equation}\label{global-version}
H^0(A,\F)\otimes H^0(A,\G)\otimes \OO_A\cong 
{d_{23}}_*(\F\boxtimes\OO_A\boxtimes\G)\rightarrow {d_{23}}_*((\F\boxtimes\OO_A\boxtimes\G)_{|\Delta_{12}})\cong \F\hat *\G
\end{equation}
where  $d_{23}(x_1,x_2,x_3)= (x_2-x_3)$ and $\Delta_{12}=\{(x_1,x_2,x_3)\>|\>x_1=x_2\}$. 
(Note that $d_*(\OO_A\boxtimes \G)$ is trivial and canonically isomorphic to $H^0(A,\G)\otimes \OO_A$, as  is most easily seen via the automorphism of $A\times A$,  $(x,y)\mapsto (x,x-y)$, sending $p_2$ to $d$ and leaving $p_1$ unchanged. Therefore ${d_{23}}_*(\F\boxtimes\OO_A\boxtimes\G)\cong H^0(A,\F)\otimes H^0(A,\G)\otimes \OO_A$.) 

More generally, given another IT(0) sheaf on $A$, say $\H$, the multiplication map of global sections
 \[
H^0(A, \F)\otimes H^0(A, \H\otimes t_z^*\G)\rightarrow H^0(A, \F\otimes\H\otimes  t_z^*\G)
 \]
is naturally identified the fiber map at $z$ of the map of $\OO_A$-modules:
\[H^0(A,\F)\otimes ( \H\hat * \G)\cong 
{d_{23}}_*(\F\boxtimes\H\boxtimes\G)\rightarrow {d_{23}}_*((\F\boxtimes\H\boxtimes\G)_{|\Delta_{12}})\cong (\F\otimes\H)\hat *\G.
\]
 
After these preliminaries, we can proceed with the proof. Let us fix $y\in A$. Then diagram (\ref{diagram}) is identified to the  diagram of fiber maps  at $z\in A$ of the commutative diagram of $\OO_A$-modules, with surjective vertical maps,
  \begin{equation}\label{diagram1}\xymatrix{H^0(\W_{n-1,n})\otimes H^0( t_y^*\OO_A(n\Theta))\otimes H^0( \W_{n-1,n})\otimes \OO_A\ar[r]\ar[d]&H^0( \W_{n-1,n})\otimes (t_y^*\OO_A(n\Theta)\mathbin{\hat{ *}} \W_{n-1,n})\ar[d]\\ 
  V_n(y)\otimes H^0( \W_{n-1,n})\otimes \OO_A\ar[r]&\bigl(\W_{n-1,n}\otimes t_y^*\OO_A(n\Theta)\bigr)\mathbin{\hat {*}} \W_{n-1,n}} \ .
  \end{equation}

  After routine calculations  (summarized below),  based on well  known results of Mukai and Oprea,  one computes
  \begin{equation}\label{c1}
   c_1(\W_{a,a+b}\mathbin{\hat{ *}} \W_{b,a+b} )=(a+b)^{2g}\utheta
  \end{equation}
  and
  \begin{equation}\label{c11}
  c_1((\W_{a,a+b}\otimes \W_{b,a+b})\mathbin{\hat{ *}} \W_{a,a+b} )=(a+b)^{g+2}a^{g-1}(a+2b)^{g-1}\utheta .
  \end{equation}
  In particular, for $a=n-1$ and $b=1$ one gets
  \begin{equation}\label{c1bis}
   c_1(\W_{n-1,n}\mathbin{\hat{ *}} \OO_A(n\Theta) )=n^{2g}\utheta
  \end{equation}
  and
  \begin{equation}\label{c11bis}
  c_1((\W_{n-1,n}\otimes \OO_A(n\Theta))\mathbin{\hat{ *}}\W_{n-1,n} )=n^{g+2}(n-1)^{g-1}(n+1)^{g-1}\utheta .
  \end{equation}
  
  Assume that $\Theta$ is irreducible, and, as above, let us fix $y\in A$.  From Corollary \ref{aggiunto} and (\ref{c1bis}) it follows that the effective divisor defined by the vanishing of the determinant of the map 
  \[
  H^0( t_y^*\OO_A(n\Theta))\otimes H^0( \W_{n-1,n})\otimes \OO_A\rightarrow  (t_y^*\OO_A(n\Theta))\mathbin{\hat {*}} \W_{n-1,n}
  \]
   is the   divisor $E_{y,n}$ of (\ref{E}) (note that if $\Theta$ is irreducible then $E_{y,n}$ has no multiple components).  
  Now assume that  $\dim V_n(y)=(n)^{2g}-((n)^2-1)^g$, which corresponds precisely to the bound in Theorem \ref{main}. Then, as shown in the previous section,  the source and target of the bottom horizontal map of (\ref{diagram1}) have the same rank. The determinant of this map vanishes on an effective divisor  $D_{y,n}$, which is invariant under translation by $n$-torsion points  (since  the vector bundles $\W_{n-1,n}$ are so,  see e.g. \cite[(14)]{oprea}). Therefore, since the support of $D_{y,n}$ is contained in $E_{y,n}$, it must be equal to $E_{y,n}$.  Hence  $c_1(D_{y,n})$ should be a multiple of $n^{2g}\utheta$. On the other hand, (\ref{c11bis}) yields that 
  \begin{equation}\label{contradiction}
  c_1(D_{y,n})=n^{g+2}(n-1)^{g-1}(n+1)^{g-1}\utheta .
  \end{equation}
  This is a contradiction  as soon as $g\ge 3$. Hence, for $g\ge 3$,  if the rank of $\Theta_y(n)$ attains the maximum for some $y\in A$ then the polarization must be reducible. 
  
    Finally we  show the computation of (\ref{c1}) and (\ref{c11}). We use the Fourier-Mukai transform $\widehat\Phi_\cP:D(A)\rightarrow D(\widehat A)$
    already invoked in \S\ref{multip-section}, and also the transform in the opposite direction  \break  $\Phi_\cP:D(\widehat A)\rightarrow D(A)$, as well as their versions at the level  of Chow rings modulo numerical equivalence $\widehat\Phi_{CH}:\mathcal A(A)\rightarrow \mathcal A(\widehat A)$  and  $\Phi_{CH}:\mathcal A ( \widehat A)\rightarrow \mathcal A(A)$ (see  \cite[Proposition 1.21]{mukaiF}). By GRR they commute with the Chern character.   We use the following facts:  (a) $\Phi_\cP\circ\widehat \Phi_\cP=(-1)_A^*[-g]$; \break
  \ (b) $ch(\W_{a,b})=a^g\exp(\frac{b}{a}\utheta)$ (this follows from (\ref{fund})); \ (c) $\Phi_{CH}(\exp\frac{b}{a}\utheta)=(\frac{b}{a})^g\exp(-\frac{a}{b}\utheta)$ (\cite[\S 3.3]{oprea});   (d) If $G$ is symmetric then $\Phi_\cP(F\otimes G)=\Phi_\cP(F)\hat *\Phi_\cP(G)[g]$ and $\widehat\Phi _\cP(\F\hat *\G)=\widehat\Phi_\cP(\F)\otimes \widehat\Phi_\cP(\G)$ (this follows from \cite[(3.7)]{mukai} using that $\F\hat*\G\cong \F* (-1)^*\G)$, where $*$ is the Pontryagin product). 
  Therefore
  \begin{eqnarray*} 
  ch(\W_{a,a+b}\mathbin{\hat{ *}} \W_{b,a+b})
 & \buildrel{(d)(b)(c)}\over =&(-1)^g\Phi_{CH}\bigl((a+b)^{g}\exp(-\frac{a}{a+b}\utheta)(a+b)^{g}\exp(-\frac{b}{a+b}\utheta)\bigr)\\
  &=&(a+b)^{2g}\exp(\utheta)
  \end{eqnarray*}
  This proves (\ref{c1}).  Moreover
   \begin{eqnarray*} 
   &&ch\bigl(\,(\W_{a,a+b}\otimes \W_{b,a+b})\mathbin{\hat{ *}}\W_{a,a+b}\bigr)=\\
  & \buildrel{(a)(b)(d)}\over=&(-1)^g\Phi_{CH}\bigl(\,\widehat \Phi_{CH}\bigl(a^g\exp(\frac{a+b}{a}\utheta)b^g\exp(\frac{a+b}{b}\utheta)\bigr)\cdot \widehat\Phi_{CH} \bigl(a^g\exp(\frac{a+b}{a}\utheta)\bigr)\,\bigr)\\
 & =&(-1)^g \Phi_{CH}\bigl(\,\widehat\Phi_{CH}\bigl((ab)^g\exp(\frac{(a+b)^2}{ab}\utheta)\bigr)\cdot  \widehat\Phi_{CH} \bigl(a^g\exp(\frac{a+b}{a}\utheta)\bigr)\,\bigr)\\
& \buildrel{(c)}\over = &(-1)^g\Phi_{CH}\bigl(\,\bigl((ab)^g\frac{(a+b)^{2g}}{(ab)^g}\exp(-\frac{ab}{(a+b)^2}\utheta)\bigr)\cdot \bigl(a^g\frac{(a+b)^g}{(a)^g}\exp(-\frac{a}{a+b}\utheta)\bigr)\,\bigr)\\
& =
&\Phi_{CH}\bigl((a+b)^{3g}\exp (-\frac{a(2b+a)}{(a+b)^2}\utheta)\bigr)\\
& \buildrel{(c)}\over = &(a+b)^ga^g(2b+a)^g\exp(\frac{(a+b)^2}{a(2b+a)}\utheta),
 \end{eqnarray*}
where in the first equality we used that the vector bundles appearing in the calculation are symmetric, so that one can neglect the $(-1)_A^*$ in the formula $\Phi_\cP\circ \widehat \Phi _\cP=(-1)^*[-g]$.   This proves~(\ref{c11}). \qed
  
 \noindent\textbf{Proof for $ \mathbf{g=2}$. }  In this case the irreducibility means that $\Theta$ is a smooth  irreducible curve of genus $2$, and $A$  is its Jacobian. Assuming this,  we claim that   the isogeny $n_A$ restricted to any translate $\Theta_y$ is birational onto its image. We postpone this for the moment (see below), and we proceed with the proof.  We note that  the $n$-torsion points in $\Theta_y$ (if any)  map all to $0$, which is therefore a point of multiplicity $\Theta_y(n)$ of the curve $n_{A,y}(\Theta_y)$. The class of the curve $n_A(\Theta_y)$ is $n^2\utheta$ (indeed $n_A^*(n_A(\Theta_y))$ is the divisor 
 $ E_{y,n}$
   of (\ref{E}), whose class is $n^4\utheta$).  Hence 
  \[
  m(A,\Theta,0):=\inf_{0\in C\subset A}\Big\{\frac{\Theta\cdot C}{\mathrm{mult}_0\>C}\Big\}\le \frac{2n^2}{\Theta_y(n) },
  \] the infimum being taken over all reduced irreducible curves $C$ in $A$ passing through $0$. But $m(A,\Theta ,0)$ is the Seshadri constant of $\Theta$ at the point $0$ (actually it is constant on all points of $A$, see \cite[\S 5.1]{pag1}),  and it is known that, for irreducible principally polarized abelian surfaces $A$,  $m(A,\Theta,0)=\frac 4 3$ (as a particular case of a more general result concerning jacobians of hyperelliptic curves, see  \cite[Theorem 7]{debarre}, where Debarre attributes it to Lazarsfeld). This proves that, if $g=2$ and $\Theta$ is irreducible,  $\Theta_y(n)\le \frac{3}{2}n^2<2n^2-1=n^4-(n^2-1)^2$. This proves  the desired bound for $g=2$. \qed
  
  Finally, for the reader's convenience, we provide a proof of the previous claim, which is however a well known fact. 
  We write  $n_{A,y}:=(n_A)_{|\Theta_y}:\Theta_y\rightarrow A$. Let $x\in \Theta_y$. The   fiber of $n_{A,y}$ at $x$ is, set theoretically, a subset of  $\Theta_y$ of the form  $\{x+\eta_1,,x+\eta_2,\dots ,x+\eta_{k(x)}\}$, with $\eta_i\in A[n]$ and $\eta_1=0$. Hence $x\in \cap_{i=1}^{k(x)}\Theta_{y+\eta_i}$. In conclusion, the points $x$ such that $k(x)>1$ belong to the set of singular points of the  effective divisor $E_{y,n}$ of (\ref{E}), which is finite  since $\Theta$ is irreducible.   This proves what was claimed.

\vskip0.5truecm  In conclusion for all $g\ge 2$, and $n\ge 2$,  if there is a $y\in A$ such that the translate $\Theta_y$ contains  the maximal number of $n$-torsion points, namely $n^{2g}-(n^2-1)^g$, then $\Theta$ must be reducible. Therefore, by the decomposition theorem for p.p.a.v.'s (\cite[Theorem 4.3.1]{BL}),  $A$ splits as the polarized product of  irreducible p.p.a.v.'s $(A_i,\utheta_i)$ for $i=1,\dots ,k$, of dimension $g_i$, with $g=\sum_{i=1}^k g_i$. Furthermore  $\Theta_y=\sum_{i=1}^kp_i^*\Theta_{i,y_i}$, where $p_i$ denotes the projection $A\rightarrow A_i$, and it follows  that, for all $i$, the translates $\Theta_{i,y_i}$ contain the maximal number of $n$-torsion points. Therefore  $g_i=1$ for all $i$, since otherwise  some of $\Theta_i$'s would be reducible. It follows also that, for all $i$, $y_i\in A_i[n]$. This concludes the proof of Theorem \ref{main}.

\section*{Acklowledgements} The author thanks:  Dragos Oprea for pointing out a gap in an earlier draft of this paper, Mihnea Popa for many conversations about semihomogeneous vector bundles a long time ago, Riccardo Salvati Manni for his encouragement and many discussions and suggestions, and the referee for very accurate remarks.

\end{document}